\documentclass[12pt,reqno]{amsart}

\usepackage{amsmath}
\usepackage{latexsym}
\usepackage{amssymb}
\usepackage{booktabs}
\usepackage{amsfonts}

\usepackage{graphics,color,floatrow}

\usepackage[colorlinks=true, linkcolor=blue]{hyperref}

\makeatletter

\renewcommand{\proof}{\par\noindent{\it Proof.\ \ }}
\def\qed{\ifmmode\square\else\nolinebreak\hfill
$\Box$\fi\par\vskip12pt}

\def\ov{\overline} 
\def\l{\langle} \def\r{\rangle}

\def\mod{{\rm mod~}}

 \def\ZZ{{\mathbb Z}}

\def\calD{{\mathcal D}}

\def\calD{{\mathcal D}}

\def\ZZ{{\mathbb Z}}

\def\mod{{\rm mod~}}

\def\K{{\bf K}}

\def\Aut{{\rm Aut}}

\def\a{\alpha} \def\b{\beta}

\def\GF{{\rm GF}}

\def\B{{\rm B}} \def\W{{\rm W}}

\newtheorem{theorem}{Theorem}[section]
\newtheorem{lemma}[theorem]{Lemma}%
\newtheorem{proposition}[theorem]{Proposition}%
\newtheorem{hypothesis}[theorem]{Hypothesis}%

\newcommand \be{\begin{equation}}
\newcommand \ee{\end{equation}}
\newcommand \bes{\begin{eqnarray}}
\newcommand \ees{\end{eqnarray}}
\newcommand \bess{\begin{eqnarray*}}
\newcommand \eess{\end{eqnarray*}}

\numberwithin{equation}{section}

\begin{document}

\setlength{\leftmargini}{1.75em}
\setlength{\leftmarginii}{1.5em}

\title[Circular regular dessins of complete bipartite multigraphs]{Circular regular dessins with underlying graphs $\K_{2^e,2^f}^{(2^d)}$: the metacyclic case}

\thanks{This work was partially supported  by NSFC: 12061092 and  Yunnan Applied Basic Research Projects: 202501AT070016.}

\author{Yuntian Wu}
\address{$^{*}$ School of Mathematics\\
Yunnan Normal University\\
Kunming, Yunnan 650500\\
P. R. China}
\email{wyt2027@163.com}

\author{Wenwen Fan}
\address{$^{*}$ School of Mathematics\\
Yunnan Normal University\\
Kunming, Yunnan 650500\\
P. R. China}
\address{Yunnan Key Laboratory of Modern Analytical Mathematics and Applications\\  Kunming, Yunnan 650500\\ P. R. China}
\email{fww0871@163.com}

\date\today

\begin{abstract}

A dessin is a $2$-cell embedding of a connected bipartite graph into an orientable closed surface.
It is {\it regular} if the colour- and orientation-preserving automorphism group acts regularly on edges, and {\it circular} if every face boundary is a simple cycle.
We classify and enumerate circular regular dessins of the complete bipartite multigraph $\K_{2^e,2^f}^{(2^d)}$ ($e\ge f\ge 2$, $d\ge 1$) whose automorphism group is {\it metacyclic}.
This extends the work of Fan (circular, metacyclic) and Cai--Hu (regular, metacyclic) from simple graphs $\K_{2^e,2^f}$ to multigraphs $\K_{2^e,2^f}^{(2^d)}$.
Such dessins correspond to non-exact bicyclic factorisations $G=\l a\r\l b\r\cong\ZZ_{2^{e+d}}\ZZ_{2^{f+d}}$, where $\l a\r\cap\l b\r\cong\ZZ_{2^d}$ and $\l ba\r\cap \l a\r=\l ba\r\cap \l b\r=1$.
In general, a group admitting such a factorisation need not be metacyclic; in this paper we impose the additional condition that $G$ is a metacyclic $2$-group.
Under this assumption, we obtain explicit enumeration formulas for the number of such dessins up to isomorphism.
Together with the known classifications for simple graphs (the case $d=0$), this completes the enumeration of the metacyclic case for $2$-power complete bipartite multigraphs.


\end{abstract}

\maketitle

\section{Introduction}

A {\it dessin} $\mathcal{D}$ is an embedding $i:\Gamma\hookrightarrow \mathbb{S}$ of a connected bipartite graph $\Gamma=(\B\cup \W, E)$ (possibly with multiple edges)
into an orientable closed surface $\mathbb{S}$, where vertices in $\B$ and $\W$ are coloured black and white, respectively, such that each component of $\mathbb{S}\setminus i(\Gamma)$ is homeomorphic to an open disc.
The graph $\Gamma$ is called the {\it underlying graph}. In particular, if the underlying graph is complete bipartite, then $\calD$ is called a {\it complete dessin}.

The study of dessins originated from Cori's concept of {\it hypermaps}, i.e., embeddings of hypergraphs into surfaces (\cite{Cori}).
The field became a central topic in contemporary mathematics largely due to Grothendieck's {\it Esquisse d'un Programme},
which established that every dessin on an orientable closed surface corresponds bijectively to an algebraic curve defined over an algebraic number field $\overline{\mathbb{Q}}$ (equivalently, a compact Riemann surface),
and that the absolute Galois group $\operatorname{Gal}(\overline{\mathbb{Q}}/\mathbb{Q})$ acts naturally on such curves.
Thus, dessins connect Riemann surfaces, algebraic curves, and Galois theory into a fruitful interdisciplinary framework, see~\cite{Grothendieck, Jones-2010, Jones-Korea, siran-how} for details.

A dessin is {\it regular} if its colour- and orientation-preserving automorphism group acts regularly on the edge set.
The Galois action on general dessins is difficult to describe explicitly.
However, for regular dessins this action becomes more amenable to combinatorial analysis.
Gonzalez-Diez and Jaikin-Zapirain (\cite{G-J}) proved that the Galois group acts faithfully on regular dessins, making them a powerful tool for studying Galois actions.
Moreover, every dessin can be obtained as a quotient of some regular dessin by a subgroup of its automorphism group.
Therefore, understanding regular dessins suffices in principle to understand all dessins.
Complete regular dessins correspond to algebraic curves that can be viewed as generalisations of Fermat curves; therefore, their classification carries significant theoretical value.
This motivates the classification and enumeration of regular complete dessins, a topic of broad interdisciplinary interest.

A more restrictive notion is that of a {\it symmetric} dessin, also known as an orientably arc-regular bipartite map: a regular dessin that admits an external automorphism transposing the vertex colors.
Symmetric complete dessins correspond exactly to orientably arc-regular embeddings of complete bipartite graphs $\K_{n,n}$.
This family has been completely classified through a series of works by Du, Jones, Kwak, Nedela and \v{S}koviera, see~\cite{DJKNS-1,DJKNS-2,Jones-2007,JNS,KK,KK-2,NSZ,Jones-2010}.
Thus the symmetric case for complete bipartite graphs is now fully understood.
A natural next step is to study the larger family of {\it regular} complete dessins~(edge-transitive complete bipartite map) of $\K_{m,n}$ without requiring the colour-swapping symmetry.
This is an open problem proposed by Gareth Jones in 2007.
In recent years, the authors and collaborators have been dedicated to the study of this classification problem.
Our initial idea was to find a decomposition method analogous to the direct product decomposition of graphs to decompose complete regular dessins.
Following this approach, it is necessary to complete the classification of the following three 'basic' families of regular dessins:
\[ (m,n)=(p^e, p^f)~\mbox{with $p$ an odd prime};~(m,n)= (2^e, 2^f)~\mbox{or}~\gcd(m,n)=1.\]
For odd prime powers, Fan, Li and Qu~\cite{Fan-2018} first classified the circular (each face boundary is a simple cycle) regular dessins,
and later Hu, Nedela and Wang~\cite{Hu K} removed the circularity restriction to complete the classification for all regular dessins.
For coprime orders, Fan, Li and Qiao~\cite{Fan-2023} classified the circular case, where we proved that the automorphism group is a split metacyclic group and gave enumeration formulas per group;
however, a full classification of all (non-circular) regular dessins remains open because the number of non-isomorphic groups depends delicately on the parameters $m,n$.
For the $2$-power case $(2^e,2^f)$ with metacyclic automorphism groups, Fan~\cite{Fan-2022} classified circular regular dessins with explicit presentations and formulas,
and later Cai and Hu~\cite{Cai-Hu} extended this to a complete enumeration of all regular dessins (not necessarily circular).

In this paper, as a further step beyond \cite{Fan-2022, Cai-Hu}, we consider the generalisation to complete bipartite multigraph  $\K_{2^e,2^f}^{(2^d)}$ with edge-multiplicity $2^d$
and continue to focus on the metacyclic case for the automorphism group.
Our main results give a complete classification and enumeration of circular regular dessins of $\K_{2^e,2^f}^{(2^d)}$ whose automorphism group is metacyclic.
This extends the metacyclic classifications of Fan and of Cai--Hu from simple graphs to multigraphs.

We now state our main result.
\begin{theorem}\label{main-theorem}
Let $\Gamma = \K_{2^e,2^f}^{(2^d)}$ be a complete bipartite multigraph with $e\ge f\ge 2$ and $d\ge 1$, and let $\mathcal{D}$ be a circular regular dessin of $\Gamma$
whose automorphism group $G=\operatorname{Aut}\mathcal{D}$ is metacyclic.
\begin{itemize}
\item[(1)] If $e>f$, then $d=1$ and $G\cong\l x,y~|~ x^{2^{e+1}}=y^{2^f}=1,~x^y=x^{-1+2^{e-f+1}} \r.$ There are exactly $2^{2f-3}$ non-isomorphic circular regular dessins,
each of genus $(2^{e}-1)(2^f-1)+(2^{f-1}-2^{e-1})$.
\item[(2)] If $e=f$, then $G$ is isomorphic to one of the following:
\begin{itemize}
\item[(i)] $G_1^+(r)=\l x,y~|~x^{2^{f}}=y^{2^{f+d}}=1,~x^y=x^{1+2^r}\r$, where $2\le r\le f$;
\item[(ii)] $G_2^+(r)=\l x,y~|~x^{2^{f+d}}=y^{2^f}=1,~x^y=x^{1+2^r}\r$, where  $\max\{2,d\}\le r\leq f+d$;
\item[(iii)] $G_3^+(m,r)=\l x,y~|~x^{2^{m}}=1,~y^{2^{2f+d-m}}=x^{2^f},~x^y=x^{1+2^r}\r$, where $f<m<f+d$ and $\max\{2,m-f\}\le r<f$
\end{itemize}
The number $N_{(2^e,2^f,2^d)}$ of non-isomorphic circular regular dessins is
\[N_{(2^e,2^f,2^d)}=\begin{cases}
2^{2f-3} & \text{if } d=1 \\[4pt]
2^{2f-2} & \text{if } d\ge 2
\end{cases}\]
and each such dessin has genus $(2^{f+d-1}-1)(2^f-1)$.
\end{itemize}
\end{theorem}

\vskip0.2in
\section{Preliminaries}\label{sec:prelim}

Let $\Gamma=\K_{m,n}^{(\lambda)}$ be a complete bipartite multi-graph with bipartite sets $B$ and $W$ of size $m$ and $n$, respectively, where $\lambda$ is the edge multiplicity of $\Gamma$.
Let $\calD=(V,E,F)$ be a regular dessin with underlying graph $\Gamma$, and let $G=\Aut\calD$ be its automorphism group, which preserves the orientation of the supporting surface and the biparts $B$ and $W$.
Then $G$ acts regularly on the edges $E$ of $\calD$.
For a vertex $\a$ of $\Gamma$, let $E(\a)$ denote the set of edges which are incident with $\a$.
Then for vertices $\beta\in B$ and $\a\in W$, we have $|E(\b)|=n\lambda$ and $|E(\a)|=m\lambda$.
Since $G$ preserves the orientation of the supporting surfaces and the biparts $B$ and $W$,
the stabilisers satisfy $G_\beta\cong\ZZ_{n\lambda}$ and $G_\a\cong\ZZ_{m\lambda}$.
Let $G_\beta=\l b\r$ and $G_\a=\l a\r$.
Because $\Gamma=\K_{m,n}^{(\lambda)}$ is $G$-edge-regular, $G_\beta$ is regular on $E(\b)$ and $G_\a$ is regular on $E(\a)$. It follows that
\[G=G_\a G_\beta=\l a\r\l b\r\cong\ZZ_{m\lambda}\ZZ_{n\lambda}, \]
 so $G$ is a {\it bicyclic} group, namely, a product of two cyclic subgroups.  As the edge multiplicity of $\Gamma$ is $\lambda$,
we further obtain
\[\l a\r\cap \l b\r\cong\ZZ_\lambda.\] In this case, $G$ is called a {\it non-exact} bicyclic group of type $(m,n,\lambda)$, and the pair $(a,b)$ is called a {\it generating pair} of $G$ of type $(m,n,\lambda)$.

Conversely, given an $(m,n,\lambda)$-bicyclic group $G=\langle a\rangle\langle b\rangle\cong\mathbb{Z}_{m\lambda}\mathbb{Z}_{n\lambda}$,
one may construct a regular dessin $\mathcal{D}=(V,E,F)$ with $\Aut\mathcal{D}\cong G$ as follows:
\begin{align*}
V &= [G:\langle a\rangle] \cup [G:\langle b\rangle],\\[2pt]
E &= \{ [\langle a\rangle g,\, g,\, \langle b\rangle g] \mid g\in G \},\\[2pt]
F &= \{ C(a,b)g \mid g\in G \},
\end{align*}
where $C(a,b)$ is the boundary cycle of a face determined by the adjacency relations among edges; we refer to ~\cite{Chen-Fan-2021} for a detailed description.
The resulting dessin is denoted by $\mathcal{D}(G,a,b)$ and is uniquely determined by $G$ and the generating pair $(a,b)$. Every regular dessin can arise in this way.

The face length of $\mathcal{D}(G,a,b)$ equals $2o(ba)$, and its genus is given by
\[\mathbf{g}=1+\frac{1}{2}\bigl(|E|-|V|-|F|\bigr)=1+\frac{|G|}{2}\Bigl(1-\frac{1}{o(a)}-\frac{1}{o(b)}-\frac{1}{o(ba)}\Bigr).\]

Since $\mathcal{D}(G,a,b)$ is $G$-face-transitive, it is {\it circular} precisely when the boundary cycle $C(a,b)$ is a simple cycle.
The subgroup $\langle ba\rangle$ fixes $C(a,b)$ and acts transitively on its black and white vertices;
therefore $C(a,b)$ is a simple cycle if and only if no nontrivial power of $ba$ fixes the white vertex $\langle a\rangle$ or the black vertex $\langle b\rangle$.
This yields the following criterion.
\begin{lemma}\label{lem:circular-criterion}
{\rm (\cite{Fan-2021})} A regular dessin $\mathcal{D}(G,a,b)$ is circular if and only if
\[\langle ba\rangle\cap\langle a\rangle=1\quad\text{and}\quad\langle ba\rangle\cap\langle b\rangle=1.\]
\end{lemma}

To sum up, we obtain the following one-to-one correspondence.

\begin{proposition}\label{prop:circular-group}
If $G$ is a finite group, then the following are equivalent:
\begin{itemize}
\item[(1)]$G\cong\Aut\calD$ for some circular regular dessin $\calD$ with underlying graph $\K_{m,n}^{(\lambda)}$;
\item[(2)]$G$ has cyclic subgroups $\l a\r\cong\ZZ_{m\lambda}$ and $\l b\r\cong\ZZ_{n\lambda}$  such that
$$G=\l a\r\l b\r, ~\mbox{where}~\l a\r\cap \l b\r\cong \ZZ_{\lambda}~\mbox{and}~\l ba\r\cap \l a\r=\l ba\r\cap \l b\r=1.$$
\end{itemize}
\end{proposition}

In the classification problem of regular dessins, we are concerned with the number of non-isomorphic regular dessins. Therefore, we next give a criterion to determine whether two regular dessins are isomorphic.
\begin{lemma}\label{lem:isomorphism-dessin}
Two regular dessins $\mathcal{D}(G,a,b)$ and $\mathcal{D}(G,a',b')$ are isomorphic if and only if there exists $\sigma\in\Aut(G)$ such that $a^{\sigma}=a'$ and $b^{\sigma}=b'$.
\end{lemma}
If $\calD(G,a,b)\cong\calD(G,a',b')$, the two triples $(G,a,b)$ and $(G,a',b')$ are called {\it equivalent.}
Consequently, the number of non-isomorphic regular dessins arising from $G$ equals the number of $\Aut(G)$-equivalence classes of generating pairs.

In this paper we focus on the case where the underlying graph is the complete bipartite multigraphs $\mathbf{K}_{2^{e},2^{f}}^{(2^{d})}$,  where $e\ge f\ge 2$.
The simple graph case $d=0$ was treated in \cite{Fan-2022} for metacyclic automorphism groups; we therefore assume $d\ge 1$.
Let $\mathcal{D}$ be a circular regular dessin with $G=\Aut\mathcal{D}$. By Proposition~\ref{prop:circular-group}, $G$ satisfies the following hypothesis:

\begin{hypothesis}\label{hyp:circular}
Let $G$ be a non-exact bicyclic group of type $(2^e,2^f,2^d)$ with $e\ge f\ge 2$ and $d\ge 1$. That is,
\[G=\langle a\rangle\langle b\rangle\cong\mathbb{Z}_{2^{e+d}}\mathbb{Z}_{2^{f+d}},\quad \langle a\rangle\cap\langle b\rangle\cong\mathbb{Z}_{2^d},\]
and moreover $G$ satisfies
\[\langle ba\rangle\cap\langle a\rangle=1\quad\text{and}\quad \langle ba\rangle\cap\langle b\rangle=1.\]

\end{hypothesis}

In order to determine the structure of the group $G$ in the next section, we need to study its properties. Below we list two basic properties, the first of which shows that $G$ must be non-abelian when $e>f$.

\begin{lemma}\label{lem:abelian-case}
If $G$ satisfies Hypothesis~\ref{hyp:circular} and is abelian, then $e=f$.
\end{lemma}
\proof
Suppose that $e>f$ and $G$ is abelian. Then
\[ (ba)^{2^{f+d}} = b^{2^{f+d}}a^{2^{f+d}} = a^{2^{f+d}} \neq 1,\]
which forces $\langle ba\rangle\cap\langle a\rangle\neq 1$,  a contradiction.
\qed

For a group $G$, let $G'$ and $\Phi(G)$ denote its commutator subgroup and Frattini subgroup, respectively.
The following lemma determines the exponent of a group satisfying Hypothesis~\ref{hyp:circular},  that is, the highest order of elements in $G$. This will play a crucial role in the classification that follows.
\begin{lemma}\label{lem:exponent}
The exponent of the group $G$ which satisfies Hypothesis~\ref{hyp:circular} is $2^{e+d}$.
\end{lemma}
\proof
Define $G^{2^i}=\langle g^{2^i}\mid g\in G\rangle$. Then $G^{2^i}$ is characteristic in $G^{2^{i-1}}$ and $G^{2^{i-1}}/G^{2^i}$ has exponent $2$.
For a finite $2$-group we have $\Phi(G)=G'G^2$, and hence $G^2\subseteq\Phi(G)$. Since $\langle a^2\rangle\langle b^2\rangle\subseteq \langle a^2,b^2\rangle\subseteq G^2$,
and both $\Phi(G)$ and $\langle a^2\rangle\langle b^2\rangle$ have order $2^{e+f+d-2}$ (because $G/\Phi(G)\cong\mathbb{Z}_2^2$), we conclude that
\[
\Phi(G)=G^2=\langle a^2,b^2\rangle=\langle a^2\rangle\langle b^2\rangle.
\]
Applying the same reasoning to $\Phi(G)$, which is itself a $(2^{e-1},2^{f-1},2^{d})$-bicyclic $2$-group, and iterating $f$ times,  as $\l b^{2^f}\r=\l a^{2^e}\r\subseteq \l a^{2^f}\r$, we obtain
\[G^{2^f} = \langle a^{2^f}\rangle\langle b^{2^f}\rangle= \langle a^{2^f}\rangle \cong \mathbb{Z}_{2^{e+d-f}}.\]
For any $g\in G$, repeated squaring yields $g^{2^f}\in G^{2^f}$, and therefore $(g^{2^f})^{2^{e+d-f}} = g^{2^{e+d}} = 1$. Thus $\exp(G)\mid 2^{e+d}$. As $o(a)=2^{e+d}$, we get $\exp(G)=2^{e+d}$.
\qed

\vskip0.1in
\section{The Structure of $G$ in the Metacyclic Case}\label{sec:structure}

Let $G=\l a\r\l b\r\cong \ZZ_{2^{e+d}}\ZZ_{2^{f+d}}$ be a non-exact bicyclic group such that
\begin{equation*}
\l a\r\cap\l b\r\cong\ZZ_{2^d},\ \l ba\r\cap\l a\r=1\text{ and }\l ba\r\cap\l b\r=1,
\end{equation*}
where $e\ge f\ge 2$ and $d\geq 1$. From now on we will assume that $G$ is metacyclic.

To determine the structure of $G$, we need to introduce a result from elementary number theory.
Let $m$ be a positive integer, and let \[U_{2^m}=\{u\in \ZZ/2^m\ZZ\mid \gcd(u,2)=1\} \] be the multiplicative group consisting of all invertible elements of the ring $\ZZ/2^m\ZZ$.
Then, by~\cite[Theorem 6.10]{number}, if $m\geq 2$, we have
\[ U_{2^m}=\{\pm1+k2^r\mid 2\le r\le m,~1\le k\le 2^{m-r},~\gcd(k,2)=1\}\cong \ZZ_2\times\ZZ_{2^{m-2}}. \]
Let \[A_1=\{u\in U_{2^m}\mid u=1+k2^r\}~\mbox{and}~ A_2=\{u\in U_{2^m}\mid u=-1+k2^r\},\] where $1\le k\le 2^{m-r}$ is odd and $r\ge 2$.
Then $A_1$ and $A_2$ are both cyclic subgroups of order $2^{m-r}$, but $A_1\ne A_2$.
If $k$ is odd, then $1+k2^r$ and $1+2^r$ are both generators of $A_1$, while $-1+k2^r$ and $-1+2^r$  are both generators of $A_2$. This yields the following conclusion:

\begin{lemma}\label{r-replace}
 If $k$ is odd, then there exist odd integers $u,v$ such that
\begin{equation*}
	(1+k2^r)^{u}\equiv 1+2^r\ (\mod 2^m)~\text{ and }~(-1+k2^r)^{v}\equiv -1+2^r\ (\mod 2^m).
\end{equation*}
\end{lemma}

\vskip0.1in
We now begin to study the structure of the group $G$ satisfying Hypothesis~\ref{hyp:circular}.

\begin{lemma}\label{origin-structure}
The presentation of $G$ is one of the following:
\begin{itemize}
	\item[(1)] $G^+:=\l x,y~|~ x^{2^{m}}=1,y^{2^{e+f+d-m}}=x^{2^t},x^y=x^{1+2^r} \r$, where $\max\{m-r,f\}\le t\le m\le e+d$ and $\max\{2,2m-e-f-d\}\le r\le m\le e+d$;
	\item[(2)] $G^{-}:=\l x,y~|~ x^{2^m}=1, y^{2^{e+f+d-m}}=x^{2^t},x^y=x^{-1+2^r} \r$, where $\max\{m-1,f\}\le t\le m\le e+d$ and $\max\{2,2m-e-f-d\}\le r\le m\le e+d$.
\end{itemize}
\end{lemma}

\proof
Because $G$ is a metacyclic $2$-group, by \cite[III.11.2]{H} it can be presented as
\[ G=\l x,y ~|~ x^{2^m}=1,~y^{2^n}=x^{2^t},~x^y=x^{\theta} \r\]
with parameters satisfying
\[ 1\le t\le m,~~\theta^{2^n}\equiv 1\ (\mod{2^m}),~~\text{and}~~2^t(\theta-1)\equiv0\ (\mod{2^m}).\]
Here $\langle x\rangle \cong \mathbb{Z}_{2^{m}}$ is normal and $G/\langle x\rangle \cong \mathbb{Z}_{2^{n}}$, hence $|G|=2^{m+n}$.
In our situation $|G|=2^{e+f+d}$, so $m+n = e+f+d$.
Since $o(x)=2^m$ and $o(y)=2^{n+m-t}=2^{e+f+d-t}$, it follows from Lemma~\ref{lem:exponent} that $m \le e+d$ and $e+f+d-t \le e+d$, and so $f\le t\le m \le e+d$.

Since $\theta^{2^n} \equiv 1~(\mod{2^m})$, we have $\theta \in U_{2^m}$. As $m \ge f \ge 2$, we write $\theta = \pm 1 + k 2^r$ with $k$ odd and $2 \le r \le m$.
The highest power of $2$ dividing $\theta^{2^n}-1=(\pm1+k 2^r)^{2^n}-1$ is $2^{n+r}$, so $n+r \ge m$, that is, $r \ge m-n=2m-e-f-d$. Thus $\max\{2,2m-e-f-d\} \le r \le m \le e+d$.
The congruence $2^t(\theta-1)\equiv0\ (\mod{2^m})$ implies $t+r \ge m$ when $\theta = 1+k2^{r}$, and $t+1 \ge m$ when $\theta = -1+k2^{r}$. Together with $f\le t\le m \le e+d$, we obtain the bounds in the lemma.

Suppose first that $\theta=1+k2^r$ with $k$ odd.  By Lemma \ref{r-replace}, there exists an odd integer $u$ such that $(1+k2^r)^u\equiv 1+2^r\pmod{2^m}$.
Replacing $x$ by $x^u$ and $y$ by $y^u$, we obtain the presentation
\begin{equation*}
G=G^+:=\l x,y ~|~ x^{2^m}=1,y^{2^{e+f+d-m}}=x^{2^t}, x^y=x^{1+2^r} \r,
\end{equation*}
as given in part~(1).

Assume now that $\theta=-1+k2^r$ with $k$ odd. By Lemma \ref{r-replace}, there exists an odd integer $v$ such that $(-1+k2^r)^v\equiv -1+2^r\pmod{2^m}$.
Replacing $y$ by $y^v$, we have $x^y=x^{-1+2^r}$. Moreover, since $t \ge m-1$ and $v$ is odd, it follows that $y^{2^{e+f+d-m}} = x^{v2^t} = x^{2^t}$. Thus we obtain the presentation
\[G=G^{-}:=\l x,y ~|~ x^{2^m}=1, y^{2^{e+f+d-m}}=x^{2^t},x^y=x^{-1+2^r} \r,\]
as given in part~(2).
\qed

Each element of $G$ can be written uniquely in the standard form $y^ix^j$ with $1\le i\le 2^{e+f+d-m}$, $1\le j\le 2^m$.
Let $\theta=\pm 1+2^r$. Multiplication of standard forms is given by
\begin{align}
y^kx^{\ell}\cdot y^ix^j=y^{k+i}x^{\ell\theta^i+j},\label{eq:ba}
\end{align}
and powers in $G$ are given by
\begin{align}
(y^ix^j)^{2^s}=y^{2^si}x^{j(1+\theta^i+\theta^{2i}+\dots+\theta^{2^{s-1}i})}=y^{i2^s} x^{j(\theta^{i2^s}-1)/(\theta^i-1)}.\label{eq:xiyj_2power}
\end{align}

Since $G=\l a\r\l b\r$, it follows that  $(a,b)=(y^ix^j,y^kx^{\ell})$ for some positive integers $i,j,k,\ell$, and so $i, j, k, \ell$ satisfy the following condition.
\begin{lemma}\label{linear-trans}
Let $a=y^ix^j$ and $b=y^kx^{\ell}$. Then
\[
\left|\begin{matrix}
i&j\\
k&\ell
\end{matrix}\right|\not\equiv 0\pmod{2}.
\]
\end{lemma}

\proof
Let $\ov{G}=G/\Phi(G)$ and denote by $\ov{g}$ the image of $g\in G$ in $\ov{G}$. Then
\[
\ov{a}=\ov{y}^i\ov{x}^j,\qquad \ov{b}=\ov{y}^k\ov{x}^\ell .
\]
Since $\langle \ov{a},\ov{b}\rangle=\langle \ov{x},\ov{y}\rangle=\ov{G}\cong\ZZ_2\times\ZZ_2$, we view $\ov{G}$ as a $2$-dimensional vector space over $\GF(2)$.
Then both $\{\ov{a},\ov{b}\}$ and $\{\ov{x},\ov{y}\}$ are bases, so the transition matrix
\[
\begin{bmatrix}
i & j\\
k & \ell
\end{bmatrix}
\]
is non-singular over $\GF(2)$. Hence its determinant is non-zero modulo $2$.
\qed

\vskip0.2in
\subsection{Two Technical Formulas}~~~~~~~~~~~~~~~~~~~~~~~~~~~~~~~~~~~
\vskip0.1in

In order to determine the concrete structure of $G$, we need to determine the value of $m$ in Lemma \ref{origin-structure},  and consequently compute powers of elements in $G$.

The formulas for computing powers of elements in $G$ have been given in~\cite[Lemma~3.3]{Fan-2022}; we will directly quote them here.
\begin{lemma}\label{G-power}
$($\cite{Fan-2022}$)$ Let $y^ix^j$ be an element of $G$ and let $s$ be a positive integer. Then
\begin{equation*}
(y^ix^j)^{2^s}=
\begin{cases}
y^{i2^s}x^{ju_i2^s} & \text{if } G=G^{+},\\
y^{i2^s}x^{jv2^s} & \text{if } G=G^{-} \text{ and } i \text{ is even},\\
y^{i2^s}x^{jw2^{r+s-1}} & \text{if } G=G^{-} \text{ and } i \text{ is odd},
\end{cases}
\end{equation*}
where $u_i$, $v$, $w$ satisfy
$u_i\equiv 1+i2^{r-1}\ (\mod{2^r}), v\equiv 1\ (\mod{2^r})$~ and $ w\equiv 1\ (\mod{2}).$
\end{lemma}

Since $G$ is metacyclic, the commutative subgroup $G'$ is cyclic.
The following lemma provides a commutativity-like property for high powers in $G$.

\begin{lemma}\label{ba-derive-compute}
If $G'\cong\ZZ_{2^{\lambda}}$, then for any $h,g\in G$ we have $(gh)^{2^{\lambda+1}}=g^{2^{\lambda+1}}h^{2^{\lambda+1}}$.
\end{lemma}

\proof  Let $h,g\in G$. Then $h=y^ix^j$ and $g=y^kx^\ell$ for some integers $i,j,k,\ell$. By Eq.~\eqref{eq:ba}, we have  $gh=y^{k+i}x^{\ell\theta^i+j},$  where $\theta=\pm1+2^r$ and $r\ge 2$.

Suppose first that $G=G^+$,  i.e.,  $\theta=1+2^r$. Then $|G'|=2^{m-r}$,  so $\lambda=m-r$.
It follows from Lemma~\ref{G-power} that $(y^ix^j)^{2^{\lambda+1}}=y^{i2^{\lambda+1}}x^{ju_i2^{\lambda+1}}$, where $u_i\equiv 1+i2^{r-1}\pmod{2^r}$.
Write $u_i=1+i2^{r-1}+q2^r$ with $q$ an integer. Then
\[(y^ix^j)^{2^{\lambda+1}}=y^{i2^{\lambda+1}}x^{j2^{\lambda+1}(1+i2^{r-1}+q2^{r})}=y^{i2^{\lambda+1}}x^{j2^{\lambda+1}}.\]
Replacing $y^i$ by $y^{k+i}$ and $x^j$ by $x^{\ell\theta^i+j}$, we obtain
\[(gh)^{2^{\lambda+1}}=(y^{k+i}x^{\ell\theta^i+j})^{2^{\lambda+1}}=y^{(k+i)2^{\lambda+1}}x^{(\ell\theta^i+j)2^{\lambda+1}}.\]
On the other hand,
\[g^{2^{\lambda+1}} h^{2^{\lambda+1}} =(y^kx^\ell)^{2^{\lambda+1}}(y^ix^j)^{2^{\lambda+1}}=( y^{k2^{\lambda+1}} x^{\ell2^{\lambda+1}})( y^{i2^{\lambda+1}} x^{j2^{\lambda+1}})
= y^{(k+i)2^{\lambda+1}} x^{(\ell\theta^{i2^{\lambda+1}} + j)2^{\lambda+1}}.\]
Since $\theta=1+2^r$ and $r\ge 2$, we have $\theta^{i}\equiv 1\pmod{2^r}$, and consequently $\theta^{i2^{\lambda+1}}\equiv \theta^i\equiv 1\pmod{2^r}$. Therefore
\[g^{2^{\lambda+1}} h^{2^{\lambda+1}} = y^{(k+i)2^{\lambda+1}} x^{(\ell + j)2^{\lambda+1}}=(gh)^{2^{\lambda+1}}.\]

Suppose now that $G=G^-$,  i.e.,  $\theta=-1+2^r$. Then $|G'|=2^{m-1}$, so $\lambda=m-1$, that is, $\lambda+1=m$.  By Lemma~\ref{G-power} we have
\[(y^ix^j)^{2^{\lambda+1}}=y^{i2^{\lambda+1}}.\]
Replacing $y^i$ by $y^{k+i}$ and $x^j$ by $x^{\ell\theta^i+j}$ gives
\[(gh)^{2^{\lambda+1}}=(y^{k+i}x^{\ell\theta^i+j})^{2^{\lambda+1}}=y^{(k+i)2^{\lambda+1}}.\]
Meanwhile,
\[g^{2^{\lambda+1}}h^{2^{\lambda+1}}=(y^kx^\ell)^{2^{\lambda+1}}(y^ix^j)^{2^{\lambda+1}}=y^{k2^{\lambda+1}}y^{i2^{\lambda+1}}=y^{(k+i)2^{\lambda+1}},\]
so the equality holds.
This completes the proof.
\qed

\vskip0.2in
\subsection{The Case $e>f$}~~~~~~~~~~~~~~~~~~~~~~~~~
\vskip0.1in
We first determine the structure of $G$ when $e>f$.

\begin{lemma}\label{e>f-m}
If $e>f$, then $m=e+d$ and $G'\cong\ZZ_{2^{e+d-1}}$.
\end{lemma}

\proof
By Lemma \ref{origin-structure}, $\l x^2\r$ is normal in $G$. Let $\ov{G}=G/\l x^2\r$, and let $\ov{x}$ and $\ov{y}$ be the images of $x$ and $y$, respectively.
Then $\ov{G}=\l \ov{x},\ov{y} \r$ with $o(\ov{x})=2$. Since $\ov{x}^{\ov{y}}=\ov{x}^{\pm1+2^r}=\ov{x}^{\pm1}=\ov{x}$, we have that $G/\l x^2\r$ is abelian.
Consequently, $ |G'|\le |\l x^2\r|=2^{m-1}\le 2^{e+d-1}.$

Suppose that $|G'|<2^{e+d-1}$. Let $(a,b)$ be a generating pair of type $(2^e,2^f,2^d)$ as in Hypothesis~\ref{hyp:circular}. By Lemma~\ref{ba-derive-compute}, we have
\begin{equation*}
(ba)^{2^{e+d-1}}=b^{2^{e+d-1}}a^{2^{e+d-1}}=a^{2^{e+d-1}}\ne 1,
\end{equation*}
which contradicts the fact that $\l ba\r\cap \l a\r=1$. Thus $|G'|=2^{e+d-1}$, and so $m=e+d$.
\qed

The following result shows that $\K_{2^e,2^f}^{(2^d)}$ with $e>f$ admits circular  regular dessins if and only if $G=G^{-}$, that is, $\theta=-1+2^r$.

\begin{lemma}\label{presentation-of-e>f}
If $e>f$, then $d=1$ and $G$ has such a presentation
\begin{equation*}
    G=\l x,y~|~ x^{2^{e+1}}=y^{2^f}=1,~x^y=x^{-1+2^{e-f+1}} \r.
\end{equation*}
\end{lemma}

\proof
Suppose $G=G^{+}$ as stated in Lemma~\ref{origin-structure}. Then $|G'|=|\langle x^{2^r}\rangle|=2^{m-r}\le 2^{e+d-2}$, contradicting Lemma~\ref{e>f-m}.
Thus $$G=G^{-}=\langle x,y \mid x^{2^{e+d}}=1,\; y^{2^{f}}=x^{2^t},\; x^{y}=x^{-1+2^r} \rangle,$$
where $e-f+d \le r\le e+d$ and $e+d-1\le t\le e+d$. In this case, we have  $f+r\ge e+d$ and $$o(y) =2^{f+e+d-t}=2^f~\mbox{or}~2^{f+1}.$$

Let $\calD=\calD(G,a,b)$ be a circular regular dessin of $\mathbf{K}_{2^e,2^f}^{(2^d)}$ with $e>f$.
By Proposition~\ref{prop:circular-group}, $G=\Aut\calD=\langle a\rangle\langle b\rangle$, where $o(a)=2^{e+d}$, $o(b)=2^{f+d}$, $\langle a\rangle\cap\langle b\rangle\cong\ZZ_{2^d}$,
and $\langle ba\rangle\cap \langle a\rangle=1=\langle ba\rangle\cap\langle b\rangle$. Write $(a,b)=(y^ix^j,y^kx^{\ell})$ for some positive integers $i,j,k,\ell$.
By Lemma~\ref{G-power},  we have
\[
1\ne a^{2^{e+d-1}}=(y^ix^j)^{2^{e+d-1}}=
\begin{cases}
y^{i2^{e+d-1}}x^{jv2^{e+d-1}}=x^{j2^{e+d-1}} & \text{if $i$ is even}\\
y^{i2^{e+d-1}}x^{jw2^{e+d+r-2}}=y^{i2^{e+d-1}} & \text{if $i$ is odd}
\end{cases}.\]
Suppose that $i$ is odd. Then $a^{2^{e+d-1}}=y^{i2^{e+d-1}}\neq 1$, which forces $o(y)=2^{e+d}>2^{f+1}$, a contradiction.
Therefore, $i$ is even and Lemma~\ref{linear-trans} then forces $j$ and $k$ to be odd.  By Lemma~\ref{G-power},  we obtain
\[1\ne b^{2^{f+d-1}}=(y^kx^{\ell})^{2^{f+d-1}}=y^{k2^{f+d-1}}x^{\ell w2^{f+d+r-2}}=y^{2^{f+d-1}}x^{\ell2^{f+d+r-2}}.\]
Suppose that $d>1$. Since $f+r\geq e+d$ and $o(y)=2^f$ or $2^{f+1}$, it follows that  $y^{2^{f+d-1}}=1$  and $x^{\ell2^{f+d+r-2}}=1$, whence $b^{2^{f+d-1}}=1$.
This contradiction forces $d=1$.

Because $\l ba\r\l a\r\subseteq \l ba,a\r=G$, we have $o(ba)\le 2^f$, and by Lemma~\ref{G-power},
\[1=(ba)^{2^f}=(y^{k+i}x^{\ell\theta^i+j})^{2^f}=y^{(k+i)2^f}x^{(\ell\theta^i+j)w2^{f+r-1}}=y^{2^f}x^{(\ell+j)2^{f+r-1}}.\]
If $f+r>e+1$, then
$$1\ne b^{2^{f}}=y^{2^{f}}\quad\text{and}\quad 1=(ba)^{2^f}=y^{2^f},$$
a contradiction. Thus $f+r=e+1$.

Replace $y$ by $yx^{2^{t-e}}$.  Then $x^{y}=x^{-1+2^r}$, and by Lemma~\ref{G-power}, we have
$$ (yx^{2^{t-e}})^{2^f} = y^{2^f} x^{w2^{f+r+t-e-1}} = y^{2^f} x^{2^t}= 1,$$
and therefore $G=\l x,y~|~ x^{2^{e+1}}=y^{2^f}=1,~x^y=x^{-1+2^{e-f+1}} \r\cong\ZZ_{2^{e+1}}:\ZZ_{2^f}$. This completes the proof.
\qed

\vskip0.2in
\subsection{The Case $e=f$}~~~~~~~~~~~~~~~~~~~
\vskip0.1in

We first show that when $e=f$, the complete bipartite multigraph $\K_{2^f,2^f}^{(2^d)}$ admits no circular regular dessin whose automorphism group is of type $G^-$.

\begin{lemma}\label{e=f}
Let $f\geq 2$. If $\K_{2^f,2^f}^{(2^d)}$ admits a circular regular dessin $\mathcal{D}$, then $\operatorname{Aut}\mathcal{D} \not\cong G^-$;
consequently, $\operatorname{Aut}\mathcal{D}$ must be of type $G^+$.
\end{lemma}

\proof Let $\calD$ be a circular regular dessin of $\mathbf{K}_{2^f,2^f}^{(2^d)}$. By Proposition~\ref{prop:circular-group}, $G=\Aut\calD$ is a non-exact bicyclic group of type $(2^f,2^f,2^d)$ such that
$$G=\langle a\rangle\langle b\rangle\cong\ZZ_{2^{f+d}}\ZZ_{2^{f+d}}~\mbox{with}~\langle a\rangle\cap\langle b\rangle\cong\ZZ_{2^d}~\mbox{and}~\langle ba\rangle\cap\langle a\rangle=1=\langle ba\rangle\cap\langle b\rangle.$$
Then by Lemma~\ref{origin-structure}, $G\cong G^+$ or $G^-$.

Suppose that $$G\cong G^-=\langle x,y \mid x^{2^{m}}=1,~y^{2^{2f+d-m}}=x^{2^{t}},~ x^{y}=x^{-1+2^{r}}\rangle,$$
where $\max\{2,2m-2f-d\}\le r\le m\le f+d$ and $\max\{m-1,f\}\le t\le m$.
Then $$o(y)=2^{2f+d-t}\leq 2^{f+d}.$$

Let $(a,b)=(y^ix^j,y^kx^\ell)$ for some positive integers $i,j,k,\ell$. By Lemma~\ref{G-power} we obtain
\begin{align}
1&\neq a^{2^{f+d-1}}=(y^ix^j)^{2^{f+d-1}}=\begin{cases}
y^{i2^{f+d-1}}x^{jv2^{f+d-1}}=x^{j2^{f+d-1}} & \text{if $i$ is even}\\
y^{i2^{f+d-1}}x^{jw2^{f+d+r-2}}=y^{2^{f+d-1}} & \text{if $i$ is odd}
\end{cases}\label{eq:e=f_a-}
\end{align}
and
\begin{align}
1&\neq b^{2^{f+d-1}}=(y^kx^{\ell})^{2^{f+d-1}}=\begin{cases}
y^{k2^{f+d-1}}x^{\ell v2^{f+d-1}}=x^{\ell2^{f+d-1}} & \text{if $k$ is even}\\
y^{k2^{f+d-1}}x^{\ell w2^{f+d+r-2}}=y^{2^{f+d-1}} & \text{if $k$ is odd}
\end{cases}.\label{eq:e=f_b-}
\end{align}
Since $\l ba\r\l a\r\subseteq \l ba,a\r=G$, we have $o(ba)\le2^f$. It follows from Equation~\eqref{eq:ba} and Lemma~\ref{G-power} that
\begin{equation}
1=(ba)^{2^f}=(y^{k+i}x^{\ell\theta^i+j})^{2^f}=\begin{cases}
y^{(k+i)2^f}x^{(\ell\theta^i+j)v2^f} & \text{if $k+i$ is even}\\
y^{(k+i)2^f}x^{(\ell\theta^i+j)w2^{f+r-1}} & \text{if $k+i$ is odd}
\end{cases},\label{eq:e=f_ba-}
\end{equation}
where $\theta=-1+2^r$ and $r\ge2$.

Assume that $i$ is even. Then Equation~\eqref{eq:e=f_a-} gives $o(x)=2^{f+d}=2^m$, so $m=f+d$ and $f+d-1\le t\le f+d$.
Consequently, we have $y^{2^f}=x^{2^t}=1$ or $x^{2^{f+d-1}}$, implying $o(y)\le 2^{f+1}$. Meanwhile, Lemma~\ref{linear-trans} forces both $j$ and $k$ to be odd.
Equation~\eqref{eq:e=f_b-} then yields $1\neq b^{2^{f+d-1}}=y^{2^{f+d-1}}$, which means $o(y)=2^{f+d}\ge 2^{f+1}$. Combining these, we obtain $d=1$ and $o(y)=2^{f+1}$.
However, since $k+i$ is odd, Equation~\eqref{eq:e=f_ba-} gives $1=(ba)^{2^f}=y^{2^f}$, a contradiction. Therefore $i$ cannot be even; it must be odd.

Now with $i$ odd, we have from Equation~\eqref{eq:e=f_a-} that $o(y)=2^{2f+d-t}=2^{f+d}$, which forces $t=f$ and consequently $m-1\le f$.
Next, suppose that $k$ is even. Then Equation~\eqref{eq:e=f_b-} gives $1\neq b^{2^{f+d-1}}=x^{\ell2^{f+d-1}}$. This implies that $\ell$ is odd and that $2^{f+d}=o(x)=2^{m}\leq 2^{f+1}$, hence $d=1$.
Since $k+i$ is odd, Equation~\eqref{eq:e=f_ba-} yields $1=(ba)^{2^f}=y^{2^{f}}$, again a contradiction. Thus $k$ must also be odd.

Since $x^{2^f}=x^{2^t}=y^{2^{2f+d-m}}$, and with both $i$ and $k$ odd, Lemma~\ref{G-power} gives
\begin{align*}
    a^{2^{f-1}} &= y^{i2^{f-1}}x^{jw2^{f+r-2}}
                = y^{i2^{f-1}}\bigl(x^{2^f}\bigr)^{jw2^{r-2}} \in \langle y \rangle \quad\text{and}\\
    b^{2^{f-1}} &= y^{k2^{f-1}}x^{\ell w2^{f+r-2}}
                = y^{k2^{f-1}}\bigl(x^{2^f}\bigr)^{\ell w2^{r-2}} \in \langle y \rangle.
\end{align*}
Therefore $\langle a^{2^{f-1}} \rangle = \langle b^{2^{f-1}} \rangle \le \langle a \rangle \cap \langle b \rangle$. Since $o(a^{2^{f-1}}) = 2^{d+1}$,
we have $|\langle a \rangle \cap \langle b \rangle| \ge 2^{d+1}$, contradicting the fact that $\langle a \rangle \cap \langle b \rangle \cong \mathbb{Z}_{2^d}.$

Both cases lead to contradictions, completing the proof. \qed

Therefore, it suffices to investigate the specific structure of $G^+$ in this case.

\begin{lemma}\label{lem:presentation-of-e=f}
If $e=f\geq 2$, then $G=G^+$ is isomorphic to the following.
\begin{itemize}
\item[(1)] $G_1^+(r)=\l x,y~|~x^{2^{f}}=y^{2^{f+d}}=1,~x^y=x^{1+2^r}\r$, where $2\le r\le f$.
\item[(2)] $G_2^+(r)=\l x,y~|~x^{2^{f+d}}=y^{2^f}=1,~x^y=x^{1+2^r}\r$, where  $\max\{2,d\}\le r\leq f+d$.
\item[(3)] $G_3^+(m,r)=\l x,y~|~x^{2^{m}}=1,~y^{2^{2f+d-m}}=x^{2^f},~x^y=x^{1+2^r}\r$, where $f<m<f+d$ and $\max\{2,m-f\}\le r<f$.
\end{itemize}
\end{lemma}

\proof Let $G=G^+$ as stated in Lemma~\ref{origin-structure}. Since $e=f$, the presentation of $G^+$ becomes $$G=G^+=\l x,y~|~ x^{2^{m}}=1,y^{2^{2f+d-m}}=x^{2^t},x^y=x^{1+2^r} \r,$$
where $\max\{m-r,f\}\le t\le m\le f+d$ and $\max\{2,2m-2f-d\}\le r\le m$.
We now distinguish three cases according to the value of $m$.
\vskip0.1in
\noindent\textbf{Case 1: $m=f$.}
Then $t=f$ and $y^{2^{f+d}}=1$. Hence $G$ can be written as
\[G_1^+(r):=\langle x,y \mid x^{2^{f}}=1,~ y^{2^{f+d}}=1,~ x^{y}=x^{1+2^{r}} \rangle\cong \ZZ_{2^f}{:}\ZZ_{2^{f+d}},\]
where $2\le r\le f$, as stated in part~(1).

\vskip0.1in
\noindent\textbf{Case 2: $m=f+d$.} Then $y^{2^{f}}=x^{2^t}$.  Set $y_1=yx^{2^{t-f}(2^{r-1}-1)}$. By Lemma~\ref{G-power},
\[y_1^{2^f}=(yx^{2^{t-f}(2^{r-1}-1)})^{2^f}=y^{2^f}x^{u_12^{t}(2^{r-1}-1)}=y^{2^f}x^{-2^t}=1\]
because $u_1\equiv 1+2^{r-1}\ (\mod 2^r)$ and $t+r\geq m$.
Replace $y$ by $y_1$, then $G$ can be reduced to
\[G_2^+(r):=\l x,y~|~x^{2^{f+d}}=1,~ y^{2^f}=1,~x^{y}=x^{1+2^r}\r\cong \ZZ_{2^{f+d}}{:}\ZZ_{2^f},\]
where $\max\{2,d\}\le r\le f+d$, as stated in part~(2).

\vskip0.1in
\noindent\textbf{Case 3: $f<m<f+d$.}
We first show that $t=f$. Suppose to the contrary that $t>f$.  Let $(a,b)=(y^ix^j,y^kx^\ell)$ be a generating pair of type $(2^f,2^f,2^d)$ that satisfies Hypothesis~\ref{hyp:circular}.
By Lemma~\ref{G-power}, we obtain
\[1\ne a^{2^{f+d-1}}=(y^ix^j)^{2^{f+d-1}}=y^{i2^{f+d-1}}x^{j u_i2^{f+d-1}}=y^{i2^{f+d-1}},\]
\[1\ne b^{2^{f+d-1}}=(y^kx^\ell)^{2^{f+d-1}}=y^{k2^{f+d-1}}x^{\ell u_k2^{f+d-1}}=y^{k2^{f+d-1}},\]
which forces $i$ and $k$ are both odd.  Lemma~\ref{linear-trans} then implies that  $j$ and $\ell$ have opposite parity.
Recall that $o(ba)\le 2^f$.  Using Eq.~\eqref{eq:ba} and Lemma~\ref{G-power},
\[1=(ba)^{2^f}=(y^{k+i}x^{\ell\theta^i+j})^{2^f}=y^{(k+i)2^f}x^{(\ell\theta^i+j)u_{k+i}2^f}, \] where $\theta=1+2^{r}$.
Since $(\ell\theta^i+j)u_{k+i}$ is odd, $m>f$ and $t>f$, it follows that $x^{(\ell\theta^i+j)u_{k+i}2^f}\neq 1$ and $x^{(\ell\theta^i+j)u_{k+i}2^f}\notin \l y\r$,
and thus $(ba)^{2^f}\ne 1$, a contradiction. Hence $t=f$, and consequently $o(y)=2^{f+d}$ and $m\le f+r$.

Now we distinguish two subcases according to $r$.

(1) If $f\le r$, then $$G'=\langle[x,y]\rangle=\langle x^{2^{r}}\rangle\le\langle x^{2^{f}}\rangle=\langle y^{2^{2f+d-m}}\rangle\le\langle y\rangle,$$
so $\langle y\rangle$ is normal in $G$ and $[x,y]=x^{2^r}=y^{k2^{f+d+r-m}}$ for some odd integer $k$.
Put $x_{1}=y$ and $y_{1}=y^{-2^{f+d-m}}x$. It follows from Lemma~\ref{G-power} that $$y_{1}^{2^{f}}=(y^{-2^{f+d-m}}x)^{2^{f}}=1.$$
Because $[x,y]=y^{k2^{f+d+r-m}}$, we have $x_1^{y_1}=y^x=x_1^{1-k2^{f+d+r-m}}$, and by Lemma~\ref{r-replace}, there exists an odd integer $u$ such that
$$(1-k2^{f+d+r-m})^{u}\equiv 1+2^{f+d+r-m}\pmod{2^{f+d}}.$$  Let $r_1=f+d+r-m$. Then $\max\{d, 2f+d-m\}\leq r_1\leq f+d$, and so $G$ has presentation such as
\[G=\langle x_{1},y_{1} \mid x_{1}^{2^{f+d}}=y_{1}^{2^{f}}=1,~ x_{1}^{y_{1}}=x_{1}^{1+2^{r_1}} \rangle\cong \ZZ_{2^{f+d}}:\ZZ_{2^f},\]
where $\max\{2, d\}\le\max\{d, 2f+d-m\}\leq r_1\leq f+d$, as stated in part~(2).

(2) If $r<f$, then $G$ has such presentation
\[G_3^+(m,r):=\l x,y~|~x^{2^{m}}=1,\ y^{2^{2f+d-m}}=x^{2^f},\ x^y=x^{1+2^r}\r,\]
where $\max\{2,m-f\}=\max\{2,2m-2f-d,m-f\}\leq r<f<m<f+d$, as stated in part~(3).
This completes the proof.\qed

\vskip0.3in
\section{Generating Pairs for $G$}\label{sec:pairs}

Let $\Omega$ denote the set of all ordered pairs $(a,b)\in G\times G$ satisfying
\begin{equation*}
G=\l a\r\l b\r\cong\ZZ_{2^{e+d}}\ZZ_{2^{f+d}},\quad \l a\r\cap\l b\r\cong\ZZ_{2^d}~\quad \text{and}\quad \l ba\r\cap \l a\r=\l ba\r\cap\l b\r=1.
\end{equation*}
According to Proposition~\ref{prop:circular-group}, each pair $(a,b)\in\Omega$ gives a circular regular dessin $\calD(G,a,b)$ whose underlying graph is $\mathbf{K}_{2^e,2^f}^{(2^d)}$ and whose automorphism group is $G$;
such a pair is called a generating pair of type $(2^e,2^f,2^d)$. Lemma~\ref{lem:isomorphism-dessin} tells us that the number of non-isomorphic circular regular dessins of $\mathbf{K}_{2^e,2^f}^{(2^d)}$ associated with $G$
equals the number of orbits of $\Aut(G)$ acting on $\Omega$.

In the present section we characterize $\Omega$; in the subsequent section we turn to $\Aut(G)$ and finally enumerate its orbits on $\Omega$.\\

\vskip0.2in
\subsection{The Case $e>f$}~~~~~~~~~~~~~~~~~~~~~~~~~~~~~~
\vskip0.1in

Let $\calD$ be a circular regular dessin with underlying graph $\mathbf{K}_{2^e,2^f}^{(2^d)}$, and let $G=\Aut\calD$ be metacyclic.
If $e>f$, it follows from Lemma~\ref{presentation-of-e>f} that $d=1$ and $G$ has the presentation
\begin{equation*}
    G=\langle x,y\mid x^{2^{e+1}}=y^{2^f}=1,\ x^y=x^{-1+2^{e-f+1}} \rangle.
\end{equation*}
Thus each generating pair $(a,b)$ of type $(2^e,2^f,2^d)$ associated with $G$ has the form $$(a,b)=(y^ix^j, y^kx^\ell),$$   where~$1\leq i,k\leq 2^f, 1\leq j,\ell\leq 2^{e+1}$ and $i,j,k,\ell$ satisfy
$\left|\begin{matrix}
i&j\\
k&\ell\end{matrix}\right|\not\equiv 0\pmod{2}.
$

In the following we determine $\Omega$ and its size when $e > f$.

\begin{lemma}\label{generating-of-e>f}
If $e>f\ge 2$, then the set $\Omega$ and its cardinality $|\Omega|$ are as follows:
\[\Omega=\{(y^{i} x^{j},  y^{k} x^{\ell}) \mid i \text{ is even and } j,k,\ell \text{ are all odd} \}~\mbox{and}~~|\Omega| = 2^{2e+2f-2}.\]
\end{lemma}
\proof
Let $G=\langle x,y\mid x^{2^{e+1}}=y^{2^f}=1,\ x^y=x^{-1+2^{e-f+1}} \rangle$ and let $(a,b)=(y^ix^j, y^kx^{\ell})$ be a generating pair associated with $G$, where~$1\leq i,k\leq 2^f, 1\leq j,\ell\leq 2^{e+1}$ and $i,j,k,\ell$ satisfy
$i\ell-kj\not\equiv 0\pmod{2}.$ Then $o(a)=2^{e+1},~o(b)=2^{f+1}$, $\l a\r\cap \l b\r\cong  \ZZ_2$ and $\l ba\r\cap \l a\r=\l ba\r\cap \l b\r=1$.

Since $o(a)=2^{e+1}$, it follows from Lemma~\ref{G-power} that
\[1\ne a^{2^{e}}=(y^ix^j)^{2^{e}}=
\begin{cases}
y^{i2^{e}}x^{jv2^{e}}=x^{j2^{e}} & \text{if $i$ is even}\\
y^{i2^{e}}x^{jw2^{2e-f}}=1 & \text{if $i$ is odd}
\end{cases},\]
and so $i$ must be even. By the determinant condition, $j$ and $k$ are both odd.
Similarly,  because $o(b) = 2^{f+1}$, we have
\[ 1\neq b^{2^{f}}=(y^kx^\ell)^{2^{f}}=y^{k2^{f}}x^{\ell w2^{f+e-f+1-1}}=x^{\ell 2^e},\]
so $\ell$ is also odd.
In this case, $a^{2^{e}}=x^{2^{e}}=b^{2^f}$, which gives $\langle a \rangle \cap \langle b \rangle \cong \mathbb{Z}_2$, as required.

Let $o(ba) = 2^{h}$. Since $\langle ba \rangle \langle a \rangle \subseteq G$, comparing orders yields $h \le f$.
Because $k+i$ is odd, Equation~\eqref{eq:ba} together with Lemma~\ref{G-power} gives
\[1 = (ba)^{2^{h}} = y^{(k+i)2^{h}} x^{(\ell \theta^{i} + j) w 2^{h+r-1}}
 = y^{(k+i)2^{h}} x^{(\ell \theta^{i} + j) w 2^{h+e-f}},\]
where $\theta = -1+2^{e-f+1}$ and $r = e-f+1$.
Since $\langle x \rangle \cap \langle y \rangle = 1$, we obtain $y^{(k+i)2^{h}} = 1$ and $x^{(\ell \theta^{i} + j) w 2^{h+e-f}} = 1$, which forces $h = f$. Moreover,
\[
1 \neq (ba)^{2^{f-1}} = y^{(k+i)2^{f-1}} x^{(\ell \theta^{i} + j) w 2^{f+r-2}}
 = y^{2^{f-1}} x^{(\ell \theta^{i} + j) 2^{e-1}},
\]
and therefore $\langle ba \rangle \cap \langle a \rangle = \langle ba \rangle \cap \langle b \rangle = 1$.

We now count $|\Omega|$. The parities of $i, j, k, \ell$ are determined as above: $i$ is even, and $j, k, \ell$ are odd.
Hence there are $2^{f-1}$ choices for $i$, $2^{e}$ choices for $j$, $2^{f-1}$ choices for $k$, and $2^{e}$ choices for $\ell$. Therefore
\[|\Omega| = 2^{f-1} \cdot 2^{e} \cdot 2^{f-1} \cdot 2^{e} = 2^{2e+2f-2}.\]
This completes the proof.\qed

\vskip0.2in
\subsection{The Case $e=f$}~~~~~~~~~~~~~~~~~~~~~~~~
\vskip0.1in

In this subsection we determine $|\Omega|$ when $e = f$.

Let $\calD$ be a circular regular dessin with underlying graph $\mathbf{K}_{2^f,2^f}^{(2^d)}$, and let $G=\Aut\calD$ be metacyclic.
By Lemmas~\ref{e=f} and~\ref{lem:presentation-of-e=f}, the group $G$ is isomorphic to one of $G_1^+(r)$, $G_2^+(r)$, or $G_3^+(m,r)$, as described in Lemma~\ref{lem:presentation-of-e=f}.

Recall that
\[\Omega=\{(a,b)\mid G=\langle a \rangle \langle b \rangle \cong \ZZ_{2^{f+d}}\ZZ_{2^{f+d}}, \langle a \rangle \cap \langle b \rangle \cong \ZZ_{2^d},
\langle ba \rangle \cap \langle a \rangle = \langle ba \rangle \cap \langle b \rangle = 1 \}.\]
Every such pair can be written as $(a,b)=(y^i x^j, y^k x^\ell)$ for some positive integers $i,j,k,\ell$.
We now determine $|\Omega|$ for each of the three families $G_1^+(r)$, $G_2^+(r)$, $G_3^+(m,r)$ by establishing the admissible ranges of $i,j,k,\ell$.

\begin{lemma}\label{generating-G_1}
If $G=G_1^+(r)$, then the set $\Omega$ and its cardinality $|\Omega|$ are as follows:
\[\Omega =\{(y^i x^j, y^k x^\ell)\mid  i,k \text{ odd},  j,\ell \text{ have opposite parity}, 2^d \mid (k+i)\}~\mbox{and}~~|\Omega|=2^{4f+d-2}.\]
\end{lemma}

\proof Let \[G=G_1^+(r)=\l x,y~|~x^{2^{f}}=y^{2^{f+d}}=1,~x^y=x^{1+2^r}\r,~~2\le r\le f.\]  Take $(a,b)=(y^ix^j, y^kx^\ell)\in \Omega$.
Then $1\leq i,k\leq 2^{f+d}, 1\leq j,\ell\leq 2^f$ and $i,j,k,\ell$ satisfy the determinant condition $i\ell-jk\not\equiv 0\pmod{2}.$
Since $o(a)=o(b)=2^{f+d}$, Lemma~\ref{G-power} gives
\[ 1\ne a^{2^{f+d-1}}=(y^ix^j)^{2^{f+d-1}}=y^{i2^{f+d-1}}~~\mbox{and}~~
1\ne b^{2^{f+d-1}}=(y^kx^\ell)^{2^{f+d-1}}=y^{k2^{f+d-1}},\]
hence $i$ and $k$ are odd. With $i,k$ odd, the determinant condition $i\ell-jk\not\equiv 0\pmod{2}$ reduces to $\ell-j\not\equiv 0\pmod{2}$, i.e., $j$ and $\ell$ have opposite parity.

Next we verify whether the above conditions on $i,j,k,\ell$ guarantee
$\langle a\rangle\cap\langle b\rangle\cong\mathbb{Z}_{2^d}$ and
$\langle ba\rangle\cap\langle a\rangle=\langle ba\rangle\cap\langle b\rangle=1$.

By Lemma~\ref{G-power}, computing the $2^f$-th powers gives
\[a^{2^{f}}=y^{i2^{f}}x^{ju_i2^{f}}=y^{i2^{f}}~~\mbox{and}~~~b^{2^{f}}=y^{k2^{f}}x^{\ell u_k2^{f}}=y^{k2^{f}},\]
hence $\langle a^{2^f}\rangle=\langle b^{2^f}\rangle=\langle y^{2^f}\rangle\cong\mathbb{Z}_{2^d}$.
For the $2^{f-1}$-th powers:
\[a^{2^{f-1}}=y^{i2^{f-1}}x^{ju_i2^{f-1}}=
\begin{cases}
y^{i2^{f-1}}x^{2^{f-1}} & \text{if $j$ is odd}\\
y^{i2^{f-1}} & \text{if $j$ is even}
\end{cases},\]
\[b^{2^{f-1}}=y^{k2^{f-1}}x^{\ell u_k2^{f-1}}=
\begin{cases}
y^{k2^{f-1}}x^{2^{f-1}} & \text{if $\ell$ is odd}\\
y^{k2^{f-1}} & \text{if $\ell$ is even}
\end{cases}.\]
Because $j$ and $\ell$ have opposite parity, exactly one of $a^{2^{f-1}}$, $b^{2^{f-1}}$ contains the factor $x^{2^{f-1}}$, so $\langle a^{2^{f-1}}\rangle\neq\langle b^{2^{f-1}}\rangle$;
consequently $\langle a\rangle\cap\langle b\rangle\cong\mathbb{Z}_{2^d}$.

Now we check that $\langle ba\rangle\cap\langle a\rangle=\langle ba\rangle\cap\langle b\rangle=1$.
Let $o(ba)=2^t$.  Since $\langle ba\rangle\langle a\rangle\subseteq\langle ba,a\rangle=G$,  comparing the orders yields $t\le f$.  Using Equation~\eqref{eq:ba} together with Lemma~\ref{G-power}, we obtain
\[1=(ba)^{2^t}=(y^{k+i}x^{\ell\theta^i+j})^{2^t}=y^{(k+i)2^t}x^{(\ell\theta^i+j)u_{k+i}2^{t}},\]where $u_{k+i}$ is odd.
Since $\langle x\rangle\cap\langle y\rangle=1$, it follows that \[y^{(k+i)2^t}=1\quad\mbox{and}\quad x^{(\ell\theta^i+j)u_{k+i}2^{t}}=1.\]
Because $j$ and $\ell$ have opposite parity, the integer $(\ell\theta^i+j)u_{k+i}$ is odd.  Hence, from the second equality we deduce $2^f \mid 2^t$, i.e., $t \ge f$.  Together with $t \le f$, this yields $t = f$.
Substituting $t=f$ into the first equality gives $2^d\mid k+i$. In this case, we have $$(ba)^{2^{f-1}}=y^{(k+i)2^{f-1}}x^{2^{f-1}}.$$
Thus $\l (ba)^{2^{f-1}}\r \neq \l a^{2^{f+d-1}}\r$ and $\l (ba)^{2^{f-1}}\r \neq \l b^{2^{f+d-1}}\r$,  and consequently  $\langle ba\rangle\cap\langle a\rangle=\langle ba\rangle\cap\langle b\rangle=1$.
Therefore
\[\Omega=\{(y^i x^j, y^k x^\ell)\mid i,k\text{ odd}, j,\ell\text{ have opposite parity},~\mbox{and}~2^d\mid k+i\}.\]

It remains to count $|\Omega|$, with parameters $1\le i,k\le 2^{f+d}$, $1\le j,\ell\le 2^f$.
For $(j,\ell)$ with opposite parity, $j$ can be chosen arbitrarily in $2^f$ ways and $\ell$ then has $2^{f-1}$ choices, giving $2^{2f-1}$ possibilities.
For $(i,k)$ with $i,k$ odd and $2^d\mid k+i$, there are $2^{f+d-1}$ choices for odd $i$; fixing $i$, the congruence $k\equiv -i\pmod{2^d}$ together with $k$ odd selects exactly one out of every $2^{d-1}$ odd numbers,
so among the $2^{f+d-1}$ odd $k$ there are $2^{f}$ admissible values.
Thus the number of ordered pairs $(i,k)$ is $2^{f+d-1}\cdot 2^{f}=2^{2f+d-1}$.
Multiplying the two counts yields $|\Omega| = 2^{2f-1} \cdot 2^{2f+d-1} = 2^{4f+d-2}$, as claimed.\qed

\begin{lemma}\label{generating-G_2}
If $G=G_2^+(r)$, then the set $\Omega$ and its cardinality $|\Omega|$ are as follows:
\[\Omega =\{(y^i x^j, y^k x^\ell)\mid  j,\ell \text{ odd},  i,k \text{ have opposite parity}, 2^d \mid j+\ell\}~\mbox{and}~~|\Omega|=2^{4f+d-2}.\]
\end{lemma}

\proof Let \[G=G_2^+(r)=\l x,y~|~x^{2^{f+d}}=y^{2^f}=1,~x^y=x^{1+2^r}\r,~~\max\{2,d\}\le r\leq f+d.\]  Take $(a,b)=(y^ix^j, y^kx^\ell)\in \Omega$.
Then $1\leq i,k\leq 2^f, 1\leq j,\ell\leq 2^{f+d}$ and $i,j,k,\ell$ satisfy the determinant condition $i\ell-jk\not\equiv 0\pmod{2}.$
Since $o(a)=o(b)=2^{f+d}$, Lemma~\ref{G-power} gives
\[ 1\ne a^{2^{f+d-1}}=(y^ix^j)^{2^{f+d-1}}=x^{j2^{f+d-1}}~~\mbox{and}~~
1\ne b^{2^{f+d-1}}=(y^kx^\ell)^{2^{f+d-1}}=x^{\ell2^{f+d-1}},\]
hence both $j$ and $\ell$ are odd. With $j,\ell$ odd, the determinant condition $i\ell-jk\not\equiv 0\pmod{2}$ reduces to $i$ and $k$ having opposite parity.

Following the same method as in Lemma~\ref{generating-G_1}, one can show that under the above conditions on $i,j,k,\ell$ we have $\langle a\rangle\cap\langle b\rangle=\langle x^{2^f}\rangle\cong\mathbb{Z}_{2^d}$.

Next we check that $\langle ba\rangle\cap\langle a\rangle=\langle ba\rangle\cap\langle b\rangle=1$.
Let $o(ba)=2^t$. From $\langle ba\rangle\langle a\rangle\subseteq\langle ba,a\rangle=G$ we have $t\le f$.
By Equation~\eqref{eq:ba} and Lemma~\ref{G-power},
\[1=(ba)^{2^t}=(y^{k+i}x^{\ell\theta^i+j})^{2^t}=y^{(k+i)2^t}x^{(\ell\theta^i+j)u_{k+i}2^{t}}\]
with $u_{k+i}$ odd. Since $\langle x\rangle\cap\langle y\rangle=1$,
\[y^{(k+i)2^t}=1\quad\mbox{and}\quad x^{(\ell\theta^i+j)u_{k+i}2^{t}}=1.\]
Because $k+i$ is odd and $o(y)=2^f$, we get $2^f\mid2^t$, hence $t\ge f$. Together with $t\le f$, this yields $t=f$.
Substituting $t=f$ into the second equality gives
\[(\ell\theta^i+j)2^{f}\equiv0\pmod{2^{f+d}}.\]
Because $r\ge d$, we have
\[(\ell\theta^i+j)2^{f}=[\ell(1+2^r)^i+j]2^{f}\equiv(\ell+j)2^{f}\pmod{2^{f+d}},\]
which implies $2^d\mid j+\ell$.

Therefore the elements of $\Omega$ are exactly those $(a,b)$ with $j,\ell$ odd, $i,k$ of opposite parity, and $2^d\mid j+\ell$.
Similarly to the computation of $\Omega$ in Lemma~\ref{generating-G_1}, there are $2^{2f-1}$ choices for $(i,k)$ and $2^{2f+d-1}$ choices for $(j,\ell)$.
Hence $|\Omega| = 2^{2f-1} \cdot 2^{2f+d-1} = 2^{4f+d-2}$, completing the proof.\qed

\begin{lemma}\label{generating-G_3}
If $G=G_3^+(m, r)$, then $|\Omega|=2^{4f+d-2}$, where
\[\Omega =\{(y^i x^j, y^k x^\ell)\mid  i,k~\text{odd},~j,\ell \text{ have opposite parity},~\mbox{and}~k+i\equiv -(\ell+j)2^{f+d-m}\ (\mod 2^d)\}.\]
\end{lemma}

\proof Let
$$G=G_3^+(m,r)=\l x,y~|~x^{2^{m}}=1,~y^{2^{2f+d-m}}=x^{2^f},~x^y=x^{1+2^r}\r$$ with $\max\{2,m-f\}\le r<f<m<f+d.$
Take $(a,b)=(y^ix^j, y^kx^\ell)\in\Omega$. Then $1\le i,k\le 2^{2f+d-m}$, $1\le j,\ell\le 2^m$, and the determinant condition $i\ell-jk\not\equiv0\pmod2$ holds.
Since $o(a)=o(b)=2^{f+d}$, Lemma~\ref{G-power} gives
\[1\ne a^{2^{f+d-1}}=(y^ix^j)^{2^{f+d-1}}=y^{i2^{f+d-1}}\quad\text{and}\quad
1\ne b^{2^{f+d-1}}=(y^kx^\ell)^{2^{f+d-1}}=y^{k2^{f+d-1}},\]
hence $i$ and $k$ are odd. With $i,k$ odd, the determinant condition $i\ell-jk\not\equiv0\pmod2$ reduces to $j$ and $\ell$ have opposite parity.

Now we verify that $\langle a\rangle\cap\langle b\rangle\cong\mathbb{Z}_{2^d}$.
First, raising to the power $2^f$ gives
\[a^{2^f}=y^{i2^f}x^{ju_i2^f}=y^{(i+ju_i2^{f+d-m})2^f}~\mbox{and}~~b^{2^f}=y^{k2^f}x^{\ell u_k2^f}=y^{(k+\ell u_k2^{f+d-m})2^f},\] where $u_i,u_k$ are certain odd integers.
Since $i,k$ are odd and $ju_i2^{f+d-m}$, $\ell u_k2^{f+d-m}$ are even (because $f+d>m$), both exponents $i+ju_i2^{f+d-m}$ and $k+\ell u_k2^{f+d-m}$ are odd.
Hence $\langle a^{2^f}\rangle=\langle b^{2^f}\rangle=\langle y^{2^f}\rangle\cong\mathbb{Z}_{2^d}$.
Next, raising to the power $2^{f-1}$:
\[a^{2^{f-1}}=y^{i2^{f-1}}x^{ju_i2^{f-1}}=
\begin{cases}
y^{i2^{f-1}}x^{2^{f-1}} & \text{if }j\text{ is odd}\\
y^{(i+ju_i2^{f+d-m})2^{f-1}} & \text{if }j\text{ is even}
\end{cases},\]
\[b^{2^{f-1}}=y^{k2^{f-1}}x^{\ell u_k2^{f-1}}=
\begin{cases}
y^{k2^{f-1}}x^{2^{f-1}} & \text{if }\ell\text{ is odd}\\
y^{(k+\ell u_k2^{f+d-m})2^{f-1}} & \text{if }\ell\text{ is even}
\end{cases}.\]
Because $j$ and $\ell$ have opposite parity, exactly one of $a^{2^{f-1}}$ and $b^{2^{f-1}}$ contains the factor $x^{2^{f-1}}$;
therefore $\langle a^{2^{f-1}}\rangle\neq\langle b^{2^{f-1}}\rangle$, and consequently $\langle a\rangle\cap\langle b\rangle\cong\mathbb{Z}_{2^d}$.

Now we check that $\langle ba\rangle\cap\langle a\rangle=\langle ba\rangle\cap\langle b\rangle=1$.
Let $o(ba)=2^t$. From $\langle ba\rangle\langle a\rangle\subseteq\langle ba,a\rangle=G$ we have $t\le f$.
Using Equation~\eqref{eq:ba} together with Lemma~\ref{G-power}, we obtain
\[1=(ba)^{2^t}=(y^{k+i}x^{\ell\theta^i+j})^{2^t}=y^{(k+i)2^t}x^{(\ell\theta^i+j)u_{k+i}2^t},\] where $u_{k+i}\equiv 1+(k+i)2^{r-1}\equiv 1\ (\mod{2^r})$ because $i+k$ is even.
Thus $$x^{(\ell\theta^i+j)u_{k+i}2^t}=y^{-(k+i)2^t}\in \l x\r\cap \l y\r=\l x^{2^f}\r,$$ which forces \[(\ell\theta^i+j)u_{k+i}2^t\equiv0\pmod{2^f}.\]
Because $(\ell\theta^i+j)u_{k+i}$ is odd, we get $t\ge f$.
Together with $t\le f$ we obtain $t=f$. Substituting $t=f$ gives
\[y^{-(k+i)2^f}=x^{(\ell\theta^i+j)u_{k+i}2^f}=y^{(\ell\theta^i+j)u_{k+i}\cdot2^{2f+d-m}}.\]
Since $o(y)=2^{f+d}$, we deduce $(k+i)2^f+(\ell\theta^i+j)u_{k+i}\cdot2^{2f+d-m}\equiv 0~(\mbox{mod}~2^{f+d}).$
Dividing by $2^f$ yields
\[k+i+(\ell\theta^i+j)u_{k+i}\cdot2^{f+d-m}\equiv0\pmod{2^d}.\]
 Using $r\ge m-f$, one has $u_{k+i}\cdot2^{f+d-m}\equiv 2^{f+d-m}~(\mbox{mod}~2^{d})$, and therefore
$$(\ell\theta^i+j)u_{k+i}\cdot2^{f+d-m}\equiv [\ell(1+2^r)^i+j]2^{f+d-m}\equiv (\ell+j)2^{f+d-m}~(\mbox{mod}~2^{d}).$$
Consequently \[k+i\equiv -(\ell+j)2^{f+d-m}\pmod{2^{d}},\] where $0<f+d-m<d$.
Finally, with $t=f$ we have $$(ba)^{2^{f-1}}=(y^{k+i}x^{\ell\theta^i+j})^{2^{f-1}}=y^{(k+i)2^{f-1}}x^{(\ell\theta^i+j)u_{k+i}2^{f-1}},$$
which contains the factor $x^{2^{f-1}}$ (as $(\ell\theta^i+j)u_{k+i}$ is odd).
Since $x^{2^{f-1}}\notin\langle y\rangle$, it follows that $\langle ba\rangle\cap\langle a\rangle=\langle ba\rangle\cap\langle b\rangle=1$.

Therefore the elements of $\Omega$ are exactly those $(a,b)$ with $i,k$ odd, $j,\ell$ of opposite parity, and the congruence $k+i\equiv -(\ell+j)2^{f+d-m}\pmod{2^d}$.

It remains to count $|\Omega|$. We count the number of admissible tuples $(i,j,k,\ell)$ with $1\le i,k\le 2^{2f+d-m}$ and $1\le j,\ell\le 2^m$.
For $(j,\ell)$ with opposite parity: $j$ can be any of the $2^m$ numbers, and $\ell$ must have the opposite parity, giving $2^{m-1}$ choices. Hence there are $2^m\cdot2^{m-1}=2^{2m-1}$ possibilities.
For $(i,k)$ with $i,k$ odd and $k\equiv -i-(\ell+j)2^{f+d-m}\pmod{2^d}$: the number of odd $i$ in $[1,2^{2f+d-m}]$ is $2^{2f+d-m-1}$. For a fixed $i$, the congruence determines $k$ uniquely modulo $2^d$. Since $2^{f+d-m}$ is even, the right-hand side $-i-(\ell+j)2^{f+d-m}$ is odd, so the solution class consists entirely of odd integers. In the interval $[1,2^{2f+d-m}]$, each residue class modulo $2^d$ appears exactly $2^{2f+d-m}/2^d = 2^{2f-m}$ times. Therefore, for each $i$, there are $2^{2f-m}$ admissible $k$. Hence the number of pairs $(i,k)$ is
\[2^{2f+d-m-1}\cdot 2^{2f-m}=2^{4f+d-2m-1}.\]
Multiplying the two counts gives
\[|\Omega| = 2^{2m-1}\cdot 2^{4f+d-2m-1}=2^{4f+d-2},\]
as claimed. \qed

\vskip0.2in
\section{Proof of the Main Theorem}\label{sec:main}

Recall that
\[ \Omega=\{(a,b)\mid G=\langle a\rangle\langle b\rangle\cong\ZZ_{2^{e+d}}\ZZ_{2^{f+d}},\ \langle a\rangle\cap\langle b\rangle\cong\ZZ_{2^d},\ \langle ba\rangle\cap\langle a\rangle=\langle ba\rangle\cap\langle b\rangle=1 \}. \]
By Lemma~\ref{lem:isomorphism-dessin}, the number of non-isomorphic circular regular dessins of $\K_{2^e,2^f}^{(2^d)}$ associated with $G$ coincides with the number of orbits of $\Aut(G)$ acting on $\Omega$.
In this section we determine $\Aut(G)$ and its action on $\Omega$.

\vskip0.1in
The automorphism group of a split metacyclic $2$-group was determined by Curran \cite{Curran}. We quote the following special cases relevant to our groups.

\begin{lemma}\label{automorphism-of-e>f}
If $G=\l x,y~|~ x^{2^{e+1}}=y^{2^f}=1,\ x^y=x^{-1+2^{e-f+1}} \r$ with $e>f\ge 2$, then
\begin{equation*}
	|\Aut(G)|=2^{2e+1}.
\end{equation*}
\end{lemma}

\begin{lemma}\label{automorphism-metacyclic-e=f}
If $G=\l x,y~|~ x^{2^{m}}=y^{2^n}=1,\ x^y=x^{1+2^r} \r$ with $2\le r\le m-1$, then
\[
|\Aut(G)|=\begin{cases}
2^{m+n+2r-1} & \text{if $m\le n$}\\
2^{3n+r-1} & \text{if $n\le r< m$}\\
2^{2n+2r-1} & \text{if $r<n<m$}
\end{cases}.
\]
\end{lemma}

For non-split metacyclic $2$-groups, the order of the automorphism group is known. We cite the following result from \cite{chen2017automorphismgroupnonsplitmetacyclic,Cai-Hu}.

\begin{lemma}\label{lem:nonsplit_automprhism}
If $G=\l x,y~|~ x^{2^m}=1,\ y^{2^{2f+d-m}}=x^{2^f},\ x^y=x^{1+2^r}\r$, where $2\le f<m<f+d$ and $\max\{2,m-f\}\le r<f$, then
\[
|\Aut(G)|=2^{3f+d+2r-m}.
\]
\end{lemma}

We now compute the number of non-isomorphic circular regular dessins of $\K_{2^e,2^f}^{(2^d)}$.
Let $\sigma\in\Aut(G)$ be an arbitrary automorphism and suppose that $(a,b)^\sigma=(a,b)$. Since $a$ and $b$ generate $G$, the automorphism $\sigma$ fixes every element of $G$;
hence $\sigma=1$. Thus $\Aut(G)$ acts semi-regularly on $\Omega$, and the number of non-isomorphic $G$-regular circular dessins of $\K_{2^e,2^f}^{(2^d)}$ equals $|\Omega|/|\Aut(G)|$.

Let $N_{(2^e,2^f,2^d)}$ denote the total number of non-isomorphic circular regular dessins of $\K_{2^e,2^f}^{(2^d)}$. Below we discuss $N_{(2^e,2^f,2^d)}$ separately for the cases $e>f$ and $e=f$.

\begin{lemma}\label{main-lemma-1}
If $e>f$, then $N_{(2^e,2^f,2^d)}=2^{2f-3}$.
\end{lemma}

\proof
As $e>f$, by Lemma~\ref{presentation-of-e>f} and \ref{generating-of-e>f} we have
\[G=G^-=\l x,y~|~ x^{2^{e+1}}=y^{2^f}=1,x^y=x^{-1+2^{e-f+1}} \r,\]
and $|\Omega|=2^{2e+2f-2}$. By Lemma~\ref{automorphism-of-e>f}, $|\Aut(G)|=2^{2e+1}$, hence $$N_{(2^e,2^f,2^d)}=\frac{2^{2e+2f-2}}{2^{2e+1}}=2^{2f-3}.$$
\qed

\begin{lemma}\label{main-lemma-2}
If $e=f$, then
$N_{(2^e,2^f,2^d)}=\begin{cases}
2^{2f-3}&\text{if $d=1$,}\\
2^{2f-2}&\text{if $d\ge 2$.}
\end{cases}$
\end{lemma}

\proof
Since $e=f$, it follows from Lemma~\ref{lem:presentation-of-e=f} that $G$ is isomorphic to one of the following:
\begin{itemize}
\item[(1)] $G_1^+(r)=\l x,y~|~x^{2^{f}}=y^{2^{f+d}}=1,~x^y=x^{1+2^r}\r$, where $2\le r\le f$.
\item[(2)] $G_2^+(r)=\l x,y~|~x^{2^{f+d}}=y^{2^f}=1,~x^y=x^{1+2^r}\r$, where  $\max\{2,d\}\le r\leq f+d$.
\item[(3)] $G_3^+(m,r)=\l x,y~|~x^{2^{m}}=1,~y^{2^{2f+d-m}}=x^{2^f},~x^y=x^{1+2^r}\r$, where $\max\{2,m-f\}\le r<f<m<f+d$.
\end{itemize}

If $G$ is abelian, then $G\cong G_1^+(f)\cong G_2^+(f+d)\cong \ZZ_{2^f}\times \ZZ_{2^{f+d}}$.  By~\cite[Theorem 3.2]{B-C-M} and Lemmas~\ref{generating-G_1}--\ref{generating-G_2}, we have
\[|\Aut(G)|=2^{4f+d-2}=|\Omega|, \] hence there is exactly one circular regular dessin associated with this abelian group.

We now count the circular regular dessins arising from non-abelian groups among the above three families.
For each non-abelian group $G$, the number of non-isomorphic dessins is $|\Omega|/|\Aut(G)|$.

\noindent\textbf{Case 1:} $G\cong G_1^+(r)$ with $2\le r\leq f-1$.
By Lemma~\ref{automorphism-metacyclic-e=f}, we have $|\Aut(G_r)|=2^{2f+d+2r-1}$.
Combining with Lemma~\ref{generating-G_1},  the number of non-isomorphic circular regular dessins associated with each non-abelian $G_1^+(r)$ is
$$\frac{|\Omega|}{|\Aut(G)|}=\frac{2^{4f+d-2}}{2^{2f+d+2r-1}}=2^{2f-2r-1}.$$
Summing over $r=2,3,\dots,f-1$  gives the total number of circular regular dessins from this family:
\[N(G_1^+(r)) := \sum_{r=2}^{f-1} 2^{2f-2r-1}= \frac{2(2^{2f-4}-1)}{3}.\]

\noindent\textbf{Case 2:}  $G\cong G_2^+(r)$ with $\max\{2,d\}\le r\leq f+d-1$.
By Lemma~\ref{generating-G_2}, we have $|\Omega|=2^{4f+d-2}$,  while by Lemma~\ref{automorphism-metacyclic-e=f}, the order of $\Aut(G_r)$ splits into two cases:
\[|\Aut(G_r)|=
\begin{cases}
2^{2f+2r-1} & \text{if } \max\{2,d\}\le r < f,\\[4pt]
2^{3f+r-1} & \text{if } \max\{f,d\}\le r\leq f+d-1.
\end{cases}\]

If $\max\{2,d\}\le r<f$, then the number of non-isomorphic circular regular dessins from this subfamily is
\[N_1(G_2^+(r)):=\sum_{r=\max\{2,d\}}^{f-1}\frac{|\Omega|}{|\Aut(G_r)|}=\sum_{r=\max\{2,d\}}^{f-1}2^{2f+d-2r-1}=\frac{2^{d+1}}{3}(2^{2(f-\max\{2,d\})}-1).\]

If $\max\{f,d\}\le r\leq f+d-1$, then the number of non-isomorphic circular regular dessins from this subfamily is
\[N_2(G_2^+(r)):=\sum_{r=\max\{f,d\}}^{f+d-1}\frac{|\Omega|}{|\Aut(G_r)|}=\sum_{r=\max\{f,d\}}^{f+d-1}2^{f+d-r-1}=2^{\min\{f,d\}}-1.\]

\noindent\textbf{Case 3:} $G\cong G_3^+(m,r)$ with $\max\{2,m-f\}\le r<f<m<f+d$.
From the range of $m$, we have $d\ge 2$ and $m-f\le \min\{d-1,f-1\}$. By Lemma~\ref{generating-G_3} and Lemma~\ref{lem:nonsplit_automprhism}, we obtain
$|\Omega| = 2^{4f+d-2}$ and $|\Aut(G_{m,r})| = 2^{3f+d+2r-m}$, respectively.
Thus for each non-abelian group $G_3^+(m,r)$, the number of corresponding non-isomorphic circular regular dessins is
\[\frac{|\Omega|}{|\Aut(G_3^+(m,r))|}=\frac{2^{4f+d-2}}{2^{3f+d+2r-m}}=2^{f+m-2r-2}.\]
Summing over $m = f+1$ to $\min\{f+d-1,2f-1\}$ and $r = \max\{2,\,m-f\}$ to $f-1$, we obtain the total number of non-isomorphic circular regular dessins from this family:
\[\begin{aligned}
N(G_3^+(m,r)):
&= \sum_{m=f+1}^{\min\{f+d-1,2f-1\}}\;\sum_{r=\max\{2,\,m-f\}}^{f-1} 2^{f+m-2r-2}\\
&= \sum_{m=f+1}^{\min\{f+d-1,2f-1\}} 2^{m-f}\frac{2^{2(f-\max\{2,\,m-f\})}-1}{3} \\
&=\begin{cases}
\sum\limits_{m=f+1}^{f+d-1} 2^{m-f}\frac{2^{2(f-\max\{2,\,m-f\})}-1}{3}= \frac{5\cdot 2^{2f-3} - 2^{2f-d+1} - 2^{d} + 2}{3} & \text{if } d < f,\\[4pt]
\sum\limits_{m=f+1}^{2f-1} 2^{m-f}\frac{2^{2(f-\max\{2,\,m-f\})}-1}{3}= \frac{5\cdot 2^{2f-3} - 3\cdot 2^{f} + 2}{3} & \text{if } d \geq f.
\end{cases}
\end{aligned}\]

We now determine the total number $N_{(2^e,2^f,2^d)}$ of circular regular dessins of $\K_{2^e,2^f}^{(2^d)},$ distinguishing three possibilities for $d$.

If $d=1$,  the family $G_3^+(m,r)$ does not appear (since it requires $d\ge2$).
Only the abelian group, the non-abelian groups $G_1^+(r)$  and  $G_2^+(r)$ contribute.
Hence
\begin{align*}
N_{(2^e,2^f,2^d)}&=1+N(G_1^+(r))+N_1(G_2^+(r))+N_2(G_2^+(r))\\
                 &=1+\frac{2^{2f-3}-2}{3}+\frac{2^{2}}{3}(2^{2f-4}-1)+2^{1}-1\\
                 &=2^{2f-3}.
\end{align*}

 If $2\le d< f$, then all three non-abelian families $G_1^+$, $G_2^+$, $G_3^+(d<f)$ together with the abelian group appear.
Therefore
\begin{align*}
N_{(2^e,2^f,2^d)}&=1+N(G_1^+(r))+N_1(G_2^+(r))+N_2(G_2^+(r))+N(G_3^+(m,r))\\
                            &=1+\frac{2^{2f-3}-2}{3}+\frac{2^{d+1}}{3}(2^{2f-2d}-1)+2^{d}-1+\frac{5\cdot 2^{2f-3} - 2^{2f-d+1} - 2^{d} + 2}{3}\\
                            &=2^{2f-2}.
\end{align*}

 If $2\leq f\le d$, then
\begin{align*}
N_{(2^e,2^f,2^d)}&=1+N(G_1^+(r))+N_2(G_2^+(r))+N(G_3^+(m,r))(d\geq f)\\
                            &=1+\frac{2^{2f-3}-2}{3}+2^{f}-1+ \frac{5\cdot 2^{2f-3} - 3\cdot 2^{f} + 2}{3}\\
                            &=2^{2f-2}.
\end{align*}

Combining the three cases we obtain
\[N_{(2^e,2^f,2^d)} =
\begin{cases}
2^{2f-3} & \text{if } d = 1,\\[4pt]
2^{2f-2} & \text{if } d \ge 2.
\end{cases}\]
This completes the proof.
\qed
\vskip0.1in

{\noindent\bf Proof of Theorem~\ref{main-theorem}:}
The first part of Theorem~\ref{main-theorem}(the case $e>f$) follows from Lemma~\ref{presentation-of-e>f} together with Lemma~\ref{main-lemma-1}.
The second part (the case $e=f$) follows from Lemma~\ref{lem:presentation-of-e=f} and Lemma~\ref{main-lemma-2}.

For each circular regular dessin $\mathcal{D}=(G,a,b)$ of $\Gamma = \K_{2^e,2^f}^{(2^d)}$, the genus $\mathbf{g}$ is given by
\[\mathbf{g}=1+\frac{|G|}{2}(1-\frac{1}{o(a)}-\frac{1}{o(b)}-\frac{1}{o(ba)}).\]
As established in Lemmas~\ref{generating-of-e>f}-\ref{generating-G_3}, in all cases we have:
\[o(a)=2^{e+d},~o(b)=2^{f+d}~\mbox{and }~o(ba)=2^f.\]
Substituting $|G|=2^{e+f+d}$ and the orders into the genus formula yields:
\[\mathbf{g} = 1 + \frac{1}{2}(2^{e+f+d}-2^{f}-2^{e}-2^{e+d})=(2^{e+d-1}-1)(2^f-1)+(2^{f-1}-2^{e-1}).\]
One checks directly that this expression simplifies to the formulas given in Theorem~\ref{main-theorem} for the two cases $e>f$ and $e=f$.

\vskip0.2in

\noindent{\bf Data availability}
\vskip0.1in
No data was used for the research described in the article.

\vskip0.2in

\noindent{\bf Declaration of competing interest}
\vskip0.1in
We have no conflict of interest.

\end{document}